\documentclass[letterpaper,11pt]{article}
\usepackage{amsmath,amsthm,amssymb,amsfonts}
\usepackage[utf8]{inputenc}
\usepackage[T1]{fontenc}
\usepackage{lmodern}
\usepackage[english]{babel}
\usepackage[hidelinks]{hyperref}

\newtheorem{prop}{Proposition}[section]
\newtheorem{deff}{Definition}[section]
\newtheorem{teo}{Theorem}[section]
\newtheorem{lema}{Lemma}[section]
\newtheorem{obs}{Remark}[section]
\newtheorem{question}{Question}[section]
\newtheorem{coro}{Corollary}[section]
\newtheorem{claim}{Claim}[section]
\DeclareMathOperator{\colors}{colors}
\DeclareMathOperator{\succs}{succ}
\DeclareMathOperator{\Fin}{Fin}

\title{Infinite Games and Ramsey Properties of $F_\sigma$ Ideals}
\author{JOS\'{E} DE JES\'{U}S PELAYO-G\'{O}MEZ}
\date{}

\begin{document}
\maketitle

\begin{abstract}
We study combinatorial properties of Borel ideals on countable sets. We construct a tall $F_\sigma$ ideal for which player II has a winning strategy in the cut-and-choose game, answering a question of J.~Zapletal. We also generalize the characterization of Ramsey properties by the random graph ideal to ideals associated with random colorings of finite subsets. Extending a construction from \cite{articulomaicol}, we obtain, for each $m\in\omega$, a tall $F_\sigma$ ideal $\mathcal{ED}_m$ satisfying $\omega\to(\mathcal{ED}_m^+)^2_{m+1}$ but not $\omega\to(\mathcal{ED}_m^+)^2_{m+2}$. We also construct a tall $F_{\sigma\delta}$ ideal $\mathcal S_\omega$ extending the Solecki ideal, with $\mathcal S_\omega^+\to(\omega,\mathcal S_\omega^+)^2_2$ but $\mathcal S_\omega^+\not\to(\mathcal S_\omega^+)^2_2$, and derive a Ramsey-type property of the Solecki ideal. Finally, we construct a tall $F_\sigma$ $K$-uniform ideal that is not Kat\v{e}tov-equivalent to $\mathcal{ED}_{\mathrm{fin}}$, answering a question of M.~Hru\v{s}\'{a}k.
\end{abstract}

\section*{Introduction}

This article concerns combinatorial properties of ideals on countable sets, particularly their Ramsey properties. We use standard set-theoretic notation as in \cite{barto,baumgartner,bam2,kunen,descriptiva}. An \emph{ideal} on a countably infinite set $X$ is a family $\mathcal I\subseteq\mathcal P(X)$ closed under subsets and finite unions, containing all finite subsets of $X$, with $X\notin\mathcal I$. Its dual filter is $\{X\setminus A:A\in\mathcal I\}$. A set $A\subseteq X$ is \emph{$\mathcal I$-positive} if $A\notin\mathcal I$, and $\mathcal I^+$ denotes the family of these sets. We usually take $X=\omega$, but use other countable sets when they make a construction more convenient.

We consider definable ideals as subsets of the Cantor space $2^X$, identifying sets with their characteristic functions. Most of our examples are $F_\sigma$. An ideal $\mathcal I$ is \emph{tall} if every infinite $A\subseteq X$ contains an infinite $B\in\mathcal I$. For ideals $\mathcal I$ on $X$ and $\mathcal J$ on $Y$, we write $\mathcal I\le_K\mathcal J$ if there is a function $f:Y\to X$ such that $f^{-1}[A]\in\mathcal J$ for every $A\in\mathcal I$. This is the \emph{Kat\v{e}tov order} \cite{katetov}. We say that two ideals are Kat\v{e}tov-equivalent when each is Kat\v{e}tov below the other.

\begin{deff}
Let $\mathcal I$ be an ideal on $\omega$, and let $n,k\ge1$. We write $\omega\to(\mathcal I^+)^n_k$ if every coloring $c:[\omega]^n\to k$ has an $\mathcal I$-positive homogeneous set. A set $A$ is \emph{$c$-homogeneous} if $c$ is constant on $[A]^n$.
\end{deff}

We study this property and several related partition relations. The first two sections use the same family of partitions of $\omega$. The last two sections have separate constructions, although they use some of the notation and game-theoretic observations introduced earlier.

\paragraph{Section 1 (Cut-and-choose game).}
We construct a tall $F_\sigma$ ideal $\mathcal{PC}$ for which player II has a winning strategy in the cut-and-choose game $G_1(\mathcal{PC})$. This answers a question of J.~Zapletal. For the tall Borel ideals studied before, either player I was known to win or the game was undecided; as far as we know, $\mathcal{PC}$ and the families built from it are the first tall Borel ideals for which player II is known to win. The construction also gives a partition relation for dominant sets.

\paragraph{Section 2 (Ramsey properties).}
The random graph ideal $\mathcal R$ characterizes the two-color Ramsey property: $\mathcal R\le_K\mathcal I$ if and only if $\omega\not\to(\mathcal I^+)^2_2$. We extend this characterization to more colors and to colorings of finite subsets of higher arity, using a family of ideals defined from random colorings. We then construct ideals $\mathcal{ED}_m$ such that $\omega\to(\mathcal{ED}_m^+)^2_{m+1}$ but $\omega\not\to(\mathcal{ED}_m^+)^2_{m+2}$. This extends the separation between two and three colors obtained in \cite{articulomaicol}.

\paragraph{Section 3 (The Solecki ideal and $\mathcal S_\omega$).}
Hru\v{s}\'{a}k asked whether the Solecki ideal $\mathcal S$ satisfies $\Omega\to(\mathcal S^+)^2_2$ \cite[Question 5.5]{michaeldiagrama}. We do not settle this question. We introduce a tall $F_{\sigma\delta}$ ideal $\mathcal S_\omega$ containing $\mathcal S$, compare it in the Kat\v{e}tov order with $\mathrm{conv}$ and $\mathrm{nwd}$, and show that $\mathcal S_\omega^+\to(\omega,\mathcal S_\omega^+)^2_2$ although $\Omega\not\to(\mathcal S_\omega^+)^2_2$. This also gives the weaker relation $\Omega\to(\omega,\mathcal S^+)^2_2$ for the Solecki ideal. We also compare the two ideals using a finite-selection version of the cut-and-choose game.

\paragraph{Section 4 ($K$-uniform ideals).}
We answer negatively the question of whether $\mathcal{ED}_{\mathrm{fin}}$ is, up to Kat\v{e}tov equivalence, the only tall $F_\sigma$ $K$-uniform ideal \cite[Question 5.11]{michaeldiagrama}. The example is an ideal of graphs on a disjoint union of finite complete graphs.

\section{The Cut and Choose Game}

\begin{deff}\label{def:game}
Let $\mathcal I$ be an ideal on $\omega$. In the \emph{cut-and-choose game} $G_1(\mathcal I)$, player I first partitions $\omega=A_0^0\cup A_1^0$. Player II chooses $i_0\in2$ and a point $n_0\in A_{i_0}^0$. At round $m\ge1$, player I partitions the previously chosen piece as $A_{i_{m-1}}^{m-1}=A_0^m\cup A_1^m$, and player II chooses $i_m\in2$ and $n_m\in A_{i_m}^m$. All partitions are into disjoint pieces; an empty piece is allowed, but player II must choose a nonempty piece. Player I wins if $\{n_m:m\in\omega\}\in\mathcal I$; otherwise player II wins. Points may be chosen more than once.
\end{deff}

\begin{center}
\begin{tabular}{c|c|c|c|c}
I & $\omega=A_0^0\cup A_1^0$ & & $A_{i_0}^0=A_0^1\cup A_1^1$ & $\cdots$\\
II & & $i_0\in2,\ n_0\in A_{i_0}^0$ & & $\cdots$
\end{tabular}
\end{center}

If player II chooses a piece belonging to $\mathcal I$, the final set belongs to $\mathcal I$: its remaining points lie in that piece, and only finitely many points have already been chosen. We will use this observation without ending the play early.

We also use the following transfer of strategies. If $\mathcal I\le_K\mathcal J$ and player I wins $G_1(\mathcal I)$, then player I wins $G_1(\mathcal J)$. Indeed, for a witness $f:Y\to X$, take preimages of the pieces in a simulated play on $X$, and answer a point $y$ by $f(y)$ in that play. The selected set in $Y$ is contained in the preimage of the selected set in $X$. Conversely, if player II wins $G_1(\mathcal I)$ and $\mathcal J\le_K\mathcal I$, pull back player I's partitions to a play on the base of $\mathcal I$, follow the winning strategy there, and return the images of the chosen points. If their images formed a $\mathcal J$-small set, its preimage would contain the $\mathcal I$-positive set produced by the strategy, a contradiction.

\begin{prop}\label{prop:not-tall}
If $\mathcal I$ is not tall, then player II has a winning strategy in $G_1(\mathcal I)$.
\end{prop}
\begin{proof}
Choose an infinite $A\subseteq\omega$ such that every infinite subset of $A$ is $\mathcal I$-positive. At each round, player II keeps a piece whose intersection with $A$ is infinite and chooses a point in that intersection different from all previously chosen points. The resulting set is an infinite subset of $A$.
\end{proof}

\begin{deff}[The random graph ideal]
Let $\langle\omega,E\rangle$ be the random graph. The \emph{random graph ideal} $\mathcal R$ is generated by its cliques and anticliques.
\end{deff}

We refer to \cite{cameron} for the random graph and to \cite{articulomaicol} for its ideal. Section 2 gives the extension property used to construct its higher-arity and multicolor analogues.

\begin{prop}\label{prop:first}
Player I has a winning strategy in $G_1(\mathcal R)$.
\end{prop}
\begin{proof}
Let $c:[\omega]^2\to2$ assign color $0$ to edges and color $1$ to nonedges. To define partitions at previously selected points as well, put $d(a,b)=c(\{a,b\})$ for $a\ne b$, and $d(a,a)=0$.

Player I first plays $\emptyset\cup\omega$. After player II chooses $n_k$ from the current piece $D$, player I plays
\[
D_0=\{m\in D:d(n_k,m)=0\},\qquad
D_1=\{m\in D:d(n_k,m)=1\}.
\]
For each distinct point $a$ chosen during the play, let $i(a)$ be the color of the piece kept immediately after its first occurrence. Every point chosen later remains in that piece. Thus, if $a\ne b$ and the first occurrence of $a$ precedes that of $b$, then $c(\{a,b\})=i(a)$. It follows that $H_i=\{a:i(a)=i\}$ is homogeneous of color $i$. The selected set is $H_0\cup H_1\in\mathcal R$, so player I wins.
\end{proof}

\begin{prop}\label{prop:second}
For an ideal $\mathcal I$ on $\omega$, $\omega\to(\mathcal I^+)^2_2$ if and only if $\mathcal R\nleq_K\mathcal I$.
\end{prop}
\begin{proof}
See \cite[Theorem 4.1]{articulomaicol}. Proposition~\ref{prop:random-critical} below gives the corresponding argument for more colors and higher arities.
\end{proof}

\begin{prop}\label{prop:third}
If player II has a winning strategy in $G_1(\mathcal I)$, then $\omega\to(\mathcal I^+)^2_2$.
\end{prop}
\begin{proof}
Otherwise, Proposition~\ref{prop:second} gives $\mathcal R\le_K\mathcal I$. By Proposition~\ref{prop:first} and the transfer of strategies, player I would also have a winning strategy in $G_1(\mathcal I)$, which is impossible.
\end{proof}

We now answer Zapletal's question by constructing a tall $F_\sigma$ ideal $\mathcal{PC}$ for which player II wins. The preceding propositions show one obstruction: such an ideal cannot be Kat\v{e}tov above $\mathcal R$. This rules out almost every known candidate. In the diagram of Borel ideals in \cite{michaeldiagrama}, all the ideals except possibly the Solecki ideal $\mathcal S$ are Kat\v{e}tov above $\mathcal R$, so none of them works, and neither does any ideal above one of them; whether $\mathcal S$ works is open, see Section 3. So the ideal we are looking for has to be small in the Kat\v{e}tov order, and this guided the construction below. It rests on the following family of partitions indexed by finite sequences, which is the central definition of this paper: the first two sections are built on it.

\begin{deff}\label{definition:fifth}
Fix a family $\langle A_s:s\in\omega^{<\omega}\rangle$ of subsets of $\omega$ satisfying:
\begin{itemize}
\item $A_\emptyset=\omega$, and $A_s$ is infinite for every $s$;
\item $\{A_{s^\frown n}:n\in\omega\}$ partitions $A_s$;
\item for any distinct $a,b\in\omega$, some level separates them: there are $k\in\omega$ and distinct $s,t\in\omega^k$ with $a\in A_s$ and $b\in A_t$.
\end{itemize}
An infinite set $M\subseteq\omega$ is \emph{major} if, for every $s\in\omega^{<\omega}$,
\[
|M\cap A_s|=\omega\ \Longrightarrow\
\bigl|\{n:|M\cap A_{s^\frown n}|=\omega\}\bigr|=\omega.
\]
A set is \emph{dominant} if it contains a major set. Write $\mathrm{Dominant}$ for the family of dominant sets, and put $\mathcal H=\mathcal P(\omega)\setminus\mathrm{Dominant}$.
\end{deff}

Such partitions can be constructed recursively. At level $k+1$, split every block of level $k$ into countably many infinite blocks, separating any two of $0,\ldots,k$ that still lie in a common block. Clearly, $\omega$ is major.

\begin{prop}\label{prop:fourth}
$\mathcal H$ is an ideal.
\end{prop}
\begin{proof}
It suffices to show that if a dominant set $B$ is partitioned into $C_0\cup C_1$, then one of the pieces is dominant. Fix a major $M\subseteq B$, replace $C_i$ by $C_i\cap M$, and let
\[
T=\{s\in\omega^{<\omega}:|M\cap A_s|=\omega\}.
\]
This is a tree containing the root, and each node has infinitely many immediate successors in $T$. Define $\phi_\alpha:T\to3$ for $\alpha\le\omega_1$ as follows. Initially, put $\phi_0(s)=i$ if $M\cap A_s\setminus C_i$ is finite, and put $\phi_0(s)=2$ otherwise. Since $M\cap A_s$ is infinite, the first case cannot hold for both colors.

Once a node receives color $0$ or $1$, retain that color at every later stage. At a successor stage, an uncolored node $s$ receives color $i$ if infinitely many of its immediate successors have color $i$; choose $0$ if both colors are available. If neither color is available, retain $2$. At a limit stage, keep any color assigned earlier and otherwise keep $2$. Each node changes at most once, and $T$ is countable, so the process stabilizes at a countable stage. In particular, $\phi_{\omega_1}$ is a fixed point of the successor rule.

We will use a selection observation. Suppose $T'\subseteq T$ is a subtree containing the root, every node of $T'$ has infinitely many immediate successors in $T'$, and $|C_i\cap A_s|=\omega$ for every $s\in T'$. Then $C_i$ contains a major set. To see this, enumerate all nodes outside $T'$ and enumerate the nodes of $T'$ with each occurring infinitely often. At stage $r$, for the prescribed node $s\in T'$, choose a fresh point of $C_i\cap A_s$ outside the blocks of the first $r$ nodes outside $T'$. This is possible: none of these excluded nodes is an ancestor of $s$, and finitely many excluded nodes can meet only finitely many of its infinitely many immediate successors in $T'$. The selected set has infinite intersection with every block indexed by $T'$ and finite intersection with every other block. It is therefore major. If there are only finitely many excluded nodes, avoid all of them from some stage onward.

\begin{claim}
If $\phi_{\omega_1}(\emptyset)=i\in2$, then $C_i$ is dominant.
\end{claim}
\begin{proof}
Induction on the first stage at which a node receives color $i$ shows that $|C_i\cap A_s|=\omega$. At stage $0$ this follows from the definition. At a later stage it follows from the infinitely many disjoint child blocks already known to have infinite intersection with $C_i$.

Every node finally colored $i$ has infinitely many children finally colored $i$. For a node colored at stage $0$, all its children in $T$ already have that color. For a node colored later, use the children that caused its first assignment. The subtree consisting of nodes all of whose initial segments have final color $i$ now satisfies the selection observation.
\end{proof}

\begin{claim}
If $\phi_{\omega_1}(\emptyset)=2$, then both $C_0$ and $C_1$ are dominant.
\end{claim}
\begin{proof}
A node finally labeled $2$ has infinite intersection with both pieces. By stabilization, all but finitely many of its children in $T$ also have label $2$: otherwise one of the two colors would occur infinitely often and the node would receive a color. Apply the selection observation to the subtree of nodes all of whose initial segments have label $2$, separately for $C_0$ and $C_1$.
\end{proof}

These claims prove the partition assertion. Dominance is preserved by taking supersets, no finite set is dominant, and $\omega$ is dominant. Thus $\mathcal H$ is a proper ideal containing all finite sets.
\end{proof}

\begin{deff}
A set $M\subseteq\omega$ is \emph{major below $s\in\omega^{<\omega}$} if $|M\cap A_s|=\omega$ and, for every $t$ extending $s$, an infinite intersection $M\cap A_t$ implies that $M\cap A_{t^\frown n}$ is infinite for infinitely many $n$. A set is \emph{dominant below $s$} if it contains a major set below $s$. Put $\mathcal H_s=\mathcal P(\omega)\setminus\mathrm{Dominant}_s$.
\end{deff}

\begin{prop}\label{prop:local-dominant}
For every $s\in\omega^{<\omega}$, $\mathcal H_s$ is an ideal. If $A$ is dominant below $s$, then it is dominant below $s^\frown n$ for infinitely many $n$.
\end{prop}
\begin{proof}
For the first assertion, apply the proof of Proposition~\ref{prop:fourth} to the subtree of nodes extending $s$, using a major subset of $A\cap A_s$. For the second, take a major $M\subseteq A$ below $s$. There are infinitely many $n$ with $|M\cap A_{s^\frown n}|=\omega$, and for each of them $M\cap A_{s^\frown n}$ is major below $s^\frown n$.
\end{proof}

\begin{deff}\label{definition:one}
A tree $T\subseteq\omega^{<\omega}$ is \emph{boundedly branching} if every node $s\in T$ has at most $|s|+1$ immediate successors. For $X\subseteq\omega$, let $T_X=\{s:X\cap A_s\ne\emptyset\}$. Define
\begin{align*}
\mathcal S_0&=\{X\subseteq\omega:\exists s\ (X\subseteq A_s\text{ and }\forall n\ |X\cap A_{s^\frown n}|=1)\},\\
\mathcal S_1&=\{X\subseteq\omega:T_X\text{ is boundedly branching}\},\\
\mathcal{PC}&=\langle\mathcal S_0\cup\mathcal S_1\rangle,
\end{align*}
where brackets denote closure under subsets and finite unions.
\end{deff}

The tree $T_X$ has no terminal nodes. We write $\succs_{T_X}(s)=\{n:s^\frown n\in T_X\}$ for its immediate successors, identified with their last coordinates.

\begin{teo}\label{theorem:one}
Player II has a winning strategy in $G_1(\mathcal{PC})$.
\end{teo}
\begin{proof}
First observe that $X$ is $\mathcal{PC}$-positive if, for every $k$, some $s\in T_X\cap\omega^k$ has at least $(k+1)^2$ immediate successors. Indeed, suppose $X$ is covered by $r$ members of $\mathcal S_0$ and $t$ members of $\mathcal S_1$. Below the finitely many support nodes of the selectors, each selector contributes at most one immediate successor at any node. Thus, at all sufficiently deep levels,
\[
|\succs_{T_X}(s)|\le r+t(|s|+1).
\]
For sufficiently large $k$, this is less than $(k+1)^2$. This proves the observation, and also shows that $\mathcal{PC}$ is proper, since every node of $T_\omega$ has infinitely many children.

Player II plays in stages. At stage $k$, maintain a node $s_k$ of length $k$ such that the current piece is dominant below $s_k$. Initially $s_0=\emptyset$. During the next $(k+1)^2$ rounds, whenever player I partitions the current piece, keep a part dominant below $s_k$, using Proposition~\ref{prop:local-dominant}. Choose a fresh point of the retained piece in a child block $A_{s_k^\frown n}$ different from those used at this stage. There are infinitely many such child blocks with infinite intersection with the retained piece.

At the end of the stage, choose a child $s_{k+1}$ of $s_k$ below which the retained piece is dominant. This starts the next stage. The final set has at least $(k+1)^2$ immediate successors at $s_k$ for every $k$, so it is $\mathcal{PC}$-positive by the observation.
\end{proof}

\begin{prop}\label{prop:fifth}
$\mathcal{PC}$ is a tall $F_\sigma$ ideal.
\end{prop}
\begin{proof}
Let $X\subseteq\omega$ be infinite. If some node $s\in T_X$ has infinitely many children in $T_X$, choose one point of $X$ from each of those child blocks. This infinite partial selector extends to a member of $\mathcal S_0$ and hence belongs to $\mathcal{PC}$.

Otherwise $T_X$ is finitely branching. Choose nested nodes $s_i$ of length $i+1$ with $|X\cap A_{s_i}|=\omega$, and distinct points $n_i\in X\cap A_{s_i}$. Such nodes exist because an infinite intersection is partitioned among finitely many nonempty child blocks. At level $d+1$, all points $n_i$ with $i\ge d$ lie in the single block $A_{s_d}$, while there are only $d$ earlier points. Thus every node of length $d$ in $T_{\{n_i:i\in\omega\}}$ has at most $d+1$ children. This gives an infinite member of $\mathcal S_1$ contained in $X$.

For the descriptive complexity, put
\[
\mathcal S_2=\{X\subseteq\omega:\exists s\ (X\subseteq A_s\text{ and }\forall n\ |X\cap A_{s^\frown n}|\le1)\}.
\]
These are precisely the subsets of members of $\mathcal S_0$.

\begin{claim}\label{clain:one}
$\mathcal S_2$ is compact.
\end{claim}
\begin{proof}
Suppose $X\notin\mathcal S_2$. At the first level at which points of $X$ separate, at least one occupied block contains two points and another contains one; otherwise $X$ would be a partial selector. Choose distinct $a,b$ in the first block and $c$ in the second. No member of $\mathcal S_2$ can contain all three: their common support would have to be at or above their separation, where two lie in the same child. Consequently $\{Y:a,b,c\in Y\}$ is an open neighborhood of $X$ disjoint from $\mathcal S_2$. Thus $\mathcal S_2$ is closed in $2^\omega$ and hence compact.
\end{proof}

\begin{claim}\label{clain:two}
$\mathcal S_1$ is compact.
\end{claim}
\begin{proof}
If $X\notin\mathcal S_1$, some node $s$ has at least $|s|+2$ distinct children in $T_X$. Choose one point of $X$ from each of those child blocks. Requiring these finitely many points to belong to a set gives an open neighborhood of $X$ disjoint from $\mathcal S_1$.
\end{proof}

The compact hereditary family $\mathcal S_1\cup\mathcal S_2$ generates $\mathcal{PC}$. For each $r$, the unions of $r$ members of this family form a compact set, by continuity of finite union. Since the family is hereditary, these unions include every set covered by $r$ generators. Taking the union over $r$ proves that $\mathcal{PC}$ is $F_\sigma$.
\end{proof}

\begin{deff}
For $X\in\mathcal I^+$, the restriction of $\mathcal I$ to $X$ is $\mathcal I\upharpoonright X=\{I\cap X:I\in\mathcal I\}$. For ideals $\mathcal I,\mathcal J$ on $\omega$, we write $\mathcal I^+\to(\mathcal J^+)^2_2$ if every coloring of pairs from an $\mathcal I$-positive set has a $\mathcal J$-positive homogeneous subset.
\end{deff}

For the parametrized families below, we consider nondecreasing functions $f:\omega\to\omega$ with $f(0)\ge1$ and $f(n)\to\infty$.

\begin{deff}\label{def:tc}
A tree $T\subseteq\omega^{<\omega}$ is \emph{bounded by $f$} if $|T\cap\omega^n|\le f(n)$ for every $n$. Put
\[
\mathcal S_3(f)=\{X:T_X\text{ is bounded by }f\},\qquad
\mathcal{TC}(f)=\langle\mathcal S_0\cup\mathcal S_3(f)\rangle.
\]
\end{deff}

\begin{prop}\label{prop:tc}
For every such function $f$, $\mathcal{TC}(f)$ is a tall $F_\sigma$ ideal.
\end{prop}
\begin{proof}
The family $\mathcal S_3(f)$ is compact: a failure at level $n$ is witnessed by points in $f(n)+1$ distinct blocks of that level. Together with $\mathcal S_2$, it gives a compact hereditary generating family, as in Proposition~\ref{prop:fifth}.

For tallness, let $X$ be infinite. If $T_X$ has an infinitely branching node, use a partial selector as before. Otherwise choose a branch $b\in\omega^\omega$ with $|X\cap A_{b\upharpoonright d}|=\omega$ for every $d$. Choose strictly increasing $d_i$ with $f(d_i)\ge i+2$, and distinct $x_i\in X\cap A_{b\upharpoonright d_i}$. At level $n$, the points with $d_i\ge n$ lie in one block. If there are $r>0$ earlier indices, then $f(n)\ge f(d_{r-1})\ge r+1$; if there are none, use $f(n)\ge1$. Thus $|T_{\{x_i:i\in\omega\}}\cap\omega^n|\le f(n)$, so the selected set belongs to $\mathcal S_3(f)$.

Finally, a finite union of selectors and members of $\mathcal S_3(f)$ has only finitely many children at every sufficiently deep node, whereas $T_\omega$ has infinitely many. Thus the generated ideal is proper.
\end{proof}

\begin{teo}\label{thm:tc-ramsey}
For every such function $f$, $\mathcal H^+\to(\mathcal{TC}(f)^+)^2_2$.
\end{teo}
\begin{proof}
Fix a dominant set $D$. It is $\mathcal{TC}(f)$-positive: it has infinite intersection with infinitely many first-level blocks, whereas a finite union of generators has uniformly bounded intersections outside finitely many such blocks. We show that player II wins $G_1(\mathcal{TC}(f)\upharpoonright D)$. First observe that $X$ is $\mathcal{TC}(f)$-positive if, at every level $k$, some node of $T_X$ has at least
\[
q(k)=(k+1)(f(k+1)+1)
\]
children. If $X$ were covered by $r$ selectors and $t$ members of $\mathcal S_3(f)$, then at sufficiently deep nodes $s$ of length $k$ its number of children would be at most $r+t f(k+1)$. For all large $k$, this is less than $q(k)$.

Now use the strategy of Theorem~\ref{theorem:one}, starting with the dominant set $D$ and choosing $q(k)$ points in distinct child blocks during stage $k$. Dominance is the invariant that allows the choices; the winning condition is positivity for $\mathcal{TC}(f)$. The observation proves that the strategy wins. Proposition~\ref{prop:third}, on the base set $D$, gives the required partition relation.
\end{proof}

We can also vary the local bound on branching. For the same functions $f$, let
\[
\mathcal B(f)=\{X:\forall s\in T_X\ |\succs_{T_X}(s)|\le f(|s|+1)\},\qquad
\mathcal{PC}(f)=\langle\mathcal S_0\cup\mathcal B(f)\rangle.
\]
In particular, $\mathcal{PC}=\mathcal{PC}(f)$ for $f(n)=\max\{1,n\}$. The family $\mathcal B(f)$ is compact by finite witnesses, just as $\mathcal S_1$ is. Since $\mathcal S_3(f)\subseteq\mathcal B(f)$, the ideal $\mathcal{PC}(f)$ is tall. The same branching estimate and the stages of length $q(k)$ used above show that it is proper and that player II wins $G_1(\mathcal{PC}(f))$. Thus both parametrized families consist of tall $F_\sigma$ ideals.

For these conventions, $\mathcal{TC}(f)\subsetneq\mathcal{PC}(f)$. To verify strictness, construct a tree with exactly $f(d+1)$ children at each node of length $d$, together with points whose branches stay in the tree. At each finite level retain the branch of the one point already assigned to each node, choose the remaining children arbitrarily, and assign a point to each new child. This maintains at most one assigned point per node of the current level. Let $X$ be the set of all assigned points. Then $X\in\mathcal B(f)$ and
\[
|T_X\cap\omega^n|=\prod_{j=1}^n f(j).
\]
The quotient of this expression by $f(n)$ tends to infinity. Every selector has finite intersection with $X$, since $T_X$ is finitely branching. If $F$ is the finite set covered by finitely many selectors, then removing $F$ deletes at most $|F|$ nodes at each level. For every fixed $t$, therefore,
\[
|T_{X\setminus F}\cap\omega^n|\ge\prod_{j=1}^n f(j)-|F|>t f(n)
\]
for all sufficiently large $n$. Thus no finite union of selectors and members of $\mathcal S_3(f)$ covers $X$.

On the other hand, choose a nondecreasing $g$ with
\[
g(n)\ge\max\left\{f(n)+1,\prod_{j=1}^n f(j)\right\},
\]
where the empty product is $1$. Then $\mathcal B(f)\subseteq\mathcal S_3(g)$, and hence $\mathcal{PC}(f)\subseteq\mathcal{TC}(g)$. We do not know for which functions $f$ and $g$ these ideals are Kat\v{e}tov-equivalent, within each family or between the two families, nor which properties of $f$ separate them in the Kat\v{e}tov order. We think these questions deserve further study.

\section{Ramsey Properties of Ideals}

\begin{deff}
For a coloring $c:[\omega]^n\to k$ and $A\subseteq\omega$, put
\[
\colors(c,A)=c[[A]^n].
\]
We write $\omega\to(\mathcal I^+)^n_{k,l}$ if every such coloring has an $\mathcal I$-positive set $A$ with $|\colors(c,A)|\le l$. When $l=1$, we omit the second subscript. Here $n,k\ge1$ and $1\le l\le k$.
\end{deff}

We extend the construction in \cite{articulomaicol} of a tall $F_\sigma$ ideal satisfying the two-color Ramsey property but not the three-color property. We first give a multicolor and higher-arity version of the characterization by the random graph ideal. In the random-coloring construction below, $n\ge2$ and $k\ge1$.

\begin{deff}[Star property]\label{definition:two}
A coloring $c:[\omega]^n\to k$ has the \emph{$(n,k)$-star property} if, for every finite $F\subseteq\omega$ and pairwise disjoint finite families $\mathcal A_0,\ldots,\mathcal A_{k-1}\subseteq[\omega]^{n-1}$, there is
\[
x\in\omega\setminus\left(F\cup\bigcup_{i<k}\bigcup\mathcal A_i\right)
\]
such that $c(a\cup\{x\})=i$ for every $i<k$ and $a\in\mathcal A_i$.
\end{deff}

For $n=k=2$, this is a form of the usual extension property of the random graph \cite{cameron}. The finite set $F$ allows us to require a fresh witness.

\begin{prop}\label{prop:random-exists}
For $n\ge2$ and $k\ge1$, there is a coloring $R_k^n:[\omega]^n\to k$ with the $(n,k)$-star property.
\end{prop}
\begin{proof}
Construct increasing finite initial segments $X_l$ of $\omega$, beginning with $X_0=n-1$, and compatible colorings of $[X_l]^n$. At stage $l$, for each map $h:[X_l]^{n-1}\to k$, add a new vertex $x_h$ and prescribe
\[
R_k^n(a\cup\{x_h\})=h(a)\qquad(a\in[X_l]^{n-1}).
\]
These prescriptions do not conflict, since each prescribed set has just one new vertex. Color the remaining new $n$-element sets arbitrarily. Let $X_{l+1}$ include all the vertices added at this stage. There is at least one new vertex at every stage, so the initial segments exhaust $\omega$.

Given finite requirements as in Definition~\ref{definition:two}, choose $l$ containing their supports and $F$. Extend the prescribed colors to a map on $[X_l]^{n-1}$. The vertex added for that map is a witness outside $X_l$.
\end{proof}

We call $R_k^n$ a random coloring of an $n$-uniform hypergraph. Its universality is expressed by the next proposition.

\begin{prop}\label{prop:sixth}
For every coloring $c:[\omega]^n\to k$, there is an injection $f:\omega\to\omega$ such that
\[
c(a)=R_k^n(f[a])\qquad(a\in[\omega]^n).
\]
\end{prop}
\begin{proof}
Choose the images of the first $n-1$ points distinctly. Suppose $f$ has been defined on $l$, where $l\ge n-1$. For each $i<k$, let
\[
\mathcal A_i=\{f[a]:a\in[l]^{n-1},\ c(a\cup\{l\})=i\}.
\]
These families are pairwise disjoint. Apply the star property with forbidden set $f[l]$ and define $f(l)$ to be the resulting witness. This keeps $f$ injective and preserves the colors of every $n$-element set whose largest element is $l$. Recursion completes the construction.
\end{proof}

\begin{deff}
A coloring $c:[\omega]^n\to k$ has the \emph{$(n,k)$-star property on $A\subseteq\omega$} if the condition in Definition~\ref{definition:two} holds with $F\subseteq A$, $\mathcal A_i\subseteq[A]^{n-1}$, and the witness required to belong to $A$. For finite pairwise disjoint families $\mathcal A_i\subseteq[A]^{n-1}$, write
\begin{align*}
W(A,\{\mathcal A_i:i<k\})=\{x\in A\setminus\textstyle\bigcup_{i<k}\bigcup\mathcal A_i:
\ &c(a\cup\{x\})=i\\[-2pt]
&\text{for all }i<k\text{ and }a\in\mathcal A_i\}.
\end{align*}
\end{deff}

Every set with the star property is infinite, since witnesses can be chosen outside any prescribed finite set.

\begin{lema}
If $c$ has the $(n,k)$-star property on $A$, then it also has this property on $W=W(A,\{\mathcal A_i:i<k\})$. In particular, $W$ is infinite.
\end{lema}
\begin{proof}
Let $F\subseteq W$ be finite and let $\mathcal A'_i\subseteq[W]^{n-1}$ be finite and pairwise disjoint. The old supports are disjoint from $W$, so the families $\mathcal A_i\cup\mathcal A'_i$ remain pairwise disjoint. Apply the star property on $A$ to these families and the forbidden set $F$. The resulting point lies in $W$ and meets all the new requirements.
\end{proof}

\begin{lema}\label{lemma:one}
If $c$ has the $(n,k)$-star property on $A=B\cup C$, where $B$ and $C$ are disjoint, then it has the property on $B$ or on $C$.
\end{lema}
\begin{proof}
Suppose it fails on both pieces. Choose finite forbidden sets $F_B\subseteq B$, $F_C\subseteq C$ and finite families $\mathcal B_i\subseteq[B]^{n-1}$, $\mathcal C_i\subseteq[C]^{n-1}$ witnessing the failures. Since $n\ge2$, the families $\mathcal B_i\cup\mathcal C_i$ are pairwise disjoint. A witness in $A$ for their union, outside $F_B\cup F_C$, would be a witness for the failed requirements in whichever piece contains it. This is a contradiction.
\end{proof}

We use this partition property to show that the following generated ideals are proper. From now on, their parameters satisfy $n\ge2$ and $1\le l<k$.

\begin{deff}
The ideal $\mathcal R_{k,l}^n$ is generated by
\[
\{A\subseteq\omega:|\colors(R_k^n,A)|\le l\}.
\]
When $l=1$, write $\mathcal R_k^n=\mathcal R_{k,1}^n$. Thus $\mathcal R_2^2$ is the usual random graph ideal.
\end{deff}

No generator has the star property: a set with that property realizes every color. The generating family is hereditary. If finitely many generators covered $\omega$, disjointifying this cover and applying Lemma~\ref{lemma:one} repeatedly would give a piece with the star property, a contradiction. Since every finite set is covered by its singletons, this indeed defines a proper ideal containing all finite sets.

\begin{prop}
$\mathcal R_{k,l}^n$ is a tall $F_\sigma$ ideal.
\end{prop}
\begin{proof}
The generating family is closed. If a set uses at least $l+1$ colors, choose an $n$-element subset realizing each of $l+1$ distinct colors. Every set containing their finite union still uses those colors. This is an open witness to failure of membership in the generating family. The family is also hereditary, so its finite unions form an $F_\sigma$ ideal, as in Proposition~\ref{prop:fifth}.

Given an infinite set, the infinite Ramsey theorem \cite{ramsey} gives an infinite homogeneous subset, which belongs to the generating family because $l\ge1$. Thus the ideal is tall.
\end{proof}

\begin{obs}
For parameters in the preceding ranges:
\begin{itemize}
\item if $l<l'<k$, then $\mathcal R_{k,l}^n\le_K\mathcal R_{k,l'}^n$;
\item if $k<k'$ and $l<k$, then $\mathcal R_{k',l}^n\le_K\mathcal R_{k,l}^n$.
\end{itemize}
\end{obs}

The first assertion follows from inclusion. For the second, regard $R_k^n$ as a coloring with $k'$ colors, leaving the extra colors unused, and apply the next proposition.

\begin{prop}\label{prop:random-critical}
For an ideal $\mathcal I$ on $\omega$, $\mathcal R_{k,l}^n\le_K\mathcal I$ if and only if $\omega\not\to(\mathcal I^+)^n_{k,l}$.
\end{prop}
\begin{proof}
Suppose $f:\omega\to\omega$ witnesses $\mathcal R_{k,l}^n\le_K\mathcal I$. Define $c:[\omega]^n\to k$ by $c(a)=R_k^n(f[a])$ when $|f[a]|=n$, and $c(a)=0$ otherwise. Let $A\subseteq\omega$ use at most $l$ colors under $c$.

\begin{claim}
$|\colors(R_k^n,f[A])|\le l$.
\end{claim}
\begin{proof}
Every $b\in[f[A]]^n$ has a preimage $a\in[A]^n$ on which $f$ is injective. Hence $R_k^n(b)=c(a)$, and $\colors(R_k^n,f[A])\subseteq\colors(c,A)$.
\end{proof}

Thus $f[A]\in\mathcal R_{k,l}^n$ and $A\subseteq f^{-1}[f[A]]\in\mathcal I$. This proves that $c$ witnesses failure of the partition relation.

Conversely, let $c$ witness that failure, so every set using at most $l$ colors belongs to $\mathcal I$. By Proposition~\ref{prop:sixth}, embed $c$ into $R_k^n$ by an injection $f$. The preimage of each generator of $\mathcal R_{k,l}^n$ uses at most $l$ colors under $c$ and therefore belongs to $\mathcal I$. The same holds for subsets of finite unions of generators, so $f$ witnesses the desired Kat\v{e}tov reduction.
\end{proof}

Hru\v{s}\'{a}k, Meza-Alc\'{a}ntara, Th\"ummel, and Uzc\'{a}tegui constructed in \cite[Theorem 5.2]{articulomaicol} an ideal $\widetilde{\mathcal{ED}}$ satisfying the two-color Ramsey property but not the three-color property. We extend their construction to a family $\mathcal{ED}_m$ with successive finite-color thresholds. We first recall the ideals on which it is based.

\begin{deff}\label{definition:three}
For $Y\subseteq\omega\times X$, let $(Y)_i=\{x\in X:(i,x)\in Y\}$ be its $i$-th fiber. Put $\Delta=\{(i,j)\in\omega^2:j<i\}$, and define
\begin{align*}
\mathcal{ED}&=\{Y\subseteq\omega^2:\exists k\ \forall i>k\ |(Y)_i|\le k\},\\
\widetilde{\mathcal{ED}}&=\{Y\subseteq\omega^3:\exists k\ \forall i\ 
[(i\le k\Rightarrow(Y)_i\in\mathcal{ED})\ \wedge\ (i>k\Rightarrow|(Y)_i|\le k)]\},\\
\mathcal{ED}_{\mathrm{fin}}&=\mathcal{ED}\upharpoonright\Delta
=\{Y\subseteq\Delta:\sup_i|(Y)_i|<\infty\}.
\end{align*}
\end{deff}

The two descriptions of $\mathcal{ED}_{\mathrm{fin}}$ agree because every fiber of $\Delta$ is finite.

\begin{prop}
Player I has a winning strategy in $G_1(\mathcal{ED})$.
\end{prop}
\begin{proof}
Player I begins by separating column $0$ from all other columns. If player II chooses column $0$, all later points stay in an $\mathcal{ED}$-small set. Otherwise let $(n_0,m_0)$ be the chosen point. At each subsequent round, separate from the retained tail all columns up to the first coordinate of the most recently chosen point. Again, choosing that finite union of columns ensures that the final set belongs to $\mathcal{ED}$. If player II always chooses the remaining tail, the first coordinates of the selected points strictly increase, so the final set has at most one point in each column and belongs to $\mathcal{ED}$.
\end{proof}

The same strategy, intersected with $\Delta$, works for $\mathcal{ED}_{\mathrm{fin}}$. For $\widetilde{\mathcal{ED}}$, use the strategy on the first coordinate. If player II chooses a finite union of first-coordinate fibers, isolate one fiber by finitely many further cuts and then use the strategy for $\mathcal{ED}$ inside it. The finitely many earlier points do not affect the outcome.

\begin{deff}\label{definition:four}
For $m\in\omega$, put $\mathcal A_m=\{A_s:|s|=m+1\}$ and
\[
\mathcal{ED}_m=\langle\mathcal A_m\cup\mathcal S_0\rangle.
\]
\end{deff}

Thus $\mathcal{ED}_m$ is generated by the blocks at level $m+1$ and the selectors at all nodes of the partition tree.

\begin{deff}
A function $\varphi:\mathcal P(\omega)\to[0,\infty]$ is a \emph{lower semicontinuous submeasure} if $\varphi(\emptyset)=0$, $\varphi(\{i\})<\infty$ for every $i$, and
\begin{itemize}
\item $\varphi(A)\le\varphi(B)$ whenever $A\subseteq B$;
\item $\varphi(A\cup B)\le\varphi(A)+\varphi(B)$;
\item $\varphi(A)=\lim_{n\to\infty}\varphi(A\cap n)$.
\end{itemize}
Write $\Fin(\varphi)=\{A:\varphi(A)<\infty\}$.
\end{deff}

If $\varphi(\omega)=\infty$, then $\Fin(\varphi)$ is a proper $F_\sigma$ ideal. Mazur's theorem states that every $F_\sigma$ ideal has this form \cite{mazur}.

\begin{prop}\label{prop:ed-fsigma}
For every $m\in\omega$, $\mathcal{ED}_m$ is a tall $F_\sigma$ ideal.
\end{prop}
\begin{proof}
The family
\[
\mathcal C_m=\{E:\exists s\in\omega^{m+1}\ E\subseteq A_s\}\cup\mathcal S_2
\]
is compact and hereditary. For the first part, failure to lie in a single block is witnessed by two points in distinct blocks; for the second, use Claim~\ref{clain:one}. Define $\varphi_m(E)$ to be the least number of members of $\mathcal C_m$ covering $E$, and put $\varphi_m(E)=\infty$ if there is no finite cover. The empty cover gives $\varphi_m(\emptyset)=0$. Monotonicity and subadditivity follow from coverings. Each sublevel $\{E:\varphi_m(E)\le r\}$ is compact, by continuity of finite union on $\mathcal C_m^r$ and heredity. Thus, if all finite initial segments of $E$ have value at most $r$, so does $E$. This proves lower semicontinuity, and $\Fin(\varphi_m)=\mathcal{ED}_m$.

To see that the ideal is proper, a finite collection of generators is supported in finitely many first-level blocks, except for root selectors, which contribute a uniformly bounded number of points to every other first-level block. Such a collection cannot cover $\omega$. Finally, each member of $\mathcal S_1$ meets only finitely many blocks at level $m+1$, so $\mathcal{PC}\subseteq\mathcal{ED}_m$. Tallness follows from Proposition~\ref{prop:fifth}.
\end{proof}

If $m\ge n$, every block at level $m+1$ is contained in a block at level $n+1$. Hence $\mathcal{ED}_m\subseteq\mathcal{ED}_n$ and $\mathcal{ED}_m\le_K\mathcal{ED}_n$.

\begin{deff}
Ideals $\mathcal I$ on $X$ and $\mathcal J$ on $Y$ are \emph{isomorphic}, written $\mathcal I\approx\mathcal J$, if some bijection $f:X\to Y$ satisfies $A\in\mathcal I$ if and only if $f[A]\in\mathcal J$ for every $A\subseteq X$.
\end{deff}

\begin{prop}\label{prop:ed-models}
$\mathcal{ED}_0\approx\mathcal{ED}$ and $\mathcal{ED}_1\approx\widetilde{\mathcal{ED}}$.
\end{prop}
\begin{proof}
For each $i$, choose a bijection $f_i:A_{(i)}\to\omega$. The map taking $x\in A_{(i)}$ to $(i,f_i(x))$ is a bijection onto $\omega^2$. It sends first-level blocks to columns and root selectors to graphs of functions; all deeper selectors are contained in one column. These are the generators of $\mathcal{ED}$ in both directions.

For the second assertion choose bijections $g_{i,j}:A_{(i,j)}\to\omega$ and send $x\in A_{(i,j)}$ to $(i,j,g_{i,j}(x))$. Within a fixed first-coordinate fiber the generators give $\mathcal{ED}$, while root selectors have at most one point per first-coordinate fiber. A finite union therefore has $\mathcal{ED}$-small fibers in finitely many exceptional coordinates and uniformly bounded finite fibers elsewhere. Conversely, such a set is covered by finitely many generators of the two kinds. This is precisely $\widetilde{\mathcal{ED}}$.
\end{proof}

In short, the family $\mathcal{ED}_m$ extends the two known ideals: $\mathcal{ED}_0$ is $\mathcal{ED}$ and $\mathcal{ED}_1$ is $\widetilde{\mathcal{ED}}$. We designed the family with this in mind.

\begin{prop}\label{prop:ed-negative}
$\mathcal R_{m+2}^2\le_K\mathcal{ED}_m$; equivalently, $\omega\not\to(\mathcal{ED}_m^+)^2_{m+2}$.
\end{prop}
\begin{proof}
For each $a\in\omega$, let $x_a\in\omega^\omega$ be the branch determined by $a\in A_{x_a\upharpoonright j}$ for all $j$. Distinct points determine distinct branches by Definition~\ref{definition:fifth}. Put
\[
\operatorname{split}(a,b)=\min\{j:x_a(j)\ne x_b(j)\},\qquad
c(\{a,b\})=\min\{\operatorname{split}(a,b),m+1\}.
\]
A homogeneous set of color $i\le m$ is a partial selector at a node of length $i$, and a homogeneous set of color $m+1$ is contained in one block at level $m+1$. Thus every homogeneous set belongs to $\mathcal{ED}_m$. Proposition~\ref{prop:random-critical} gives the reduction. The homogeneous sets for this coloring also generate $\mathcal{ED}_m$: each block and each selector is homogeneous, with selectors below level $m$ receiving the last color.
\end{proof}

\begin{deff}
We write $\omega\to(<\omega,\ldots,<\omega\mid\mathcal I^+)^2_k$ if, for every $c:[\omega]^2\to k$, either there is an $\mathcal I$-positive homogeneous set, or every color has finite homogeneous sets of arbitrarily large size. A set is homogeneous of color $i$ when all its pairs have color $i$; this includes sets with at most one point.
\end{deff}

\begin{lema}\label{lemma:two}
For every $n,m\in\omega$, $A_{(n)}$ is $\mathcal{ED}_{m+1}$-positive and
\[
\mathcal{ED}_{m+1}\upharpoonright A_{(n)}\approx\mathcal{ED}_m.
\]
\end{lema}
\begin{proof}
For each $s\in\omega^{m+1}$ choose a bijection $f_s:A_{(n)^\frown s}\to A_s$, and let $f:A_{(n)}\to\omega$ be their union. Blocks at the terminal levels correspond under $f$. A selector above those levels is sent to a partial selector, while a selector at or below the terminal level is contained in a terminal block. A root selector intersects $A_{(n)}$ in at most one point. These observations apply in both directions, so $E\in\mathcal{ED}_{m+1}\upharpoonright A_{(n)}$ if and only if $f[E]\in\mathcal{ED}_m$. Since the latter ideal is proper, $A_{(n)}$ is positive.
\end{proof}

In plain words, $\mathcal{ED}_{m+1}$ contains countably many disjoint copies of $\mathcal{ED}_m$, one on each first-level block. This is what makes the induction in Proposition~\ref{prop:ed-positive} possible.

\begin{lema}\label{lemma:four}
If $M$ is major and $m\in\omega$, then $M$ is $\mathcal{ED}_m$-positive and $\mathcal{ED}_m\upharpoonright M\approx\mathcal{ED}_m$.
\end{lema}
\begin{proof}
We prove the assertion by induction on $m$, simultaneously for every family of partitions as in Definition~\ref{definition:fifth} on a countably infinite base set. In particular, the partitions $\langle A_{(i)^\frown s}:s\in\omega^{<\omega}\rangle$ form such a family on $A_{(i)}$. For $m>0$, the ideal $\mathcal{ED}_{m-1}$ defined from this family is exactly $\mathcal{ED}_m\upharpoonright A_{(i)}$: root selectors contribute only singletons to the restriction, and all other generators are the local generators. The bijection in Lemma~\ref{lemma:two} identifies this local ideal with $\mathcal{ED}_{m-1}$ on the original base.

First note the following description of membership. For $m=0$, $E\in\mathcal{ED}_0$ exactly when $|E\cap A_{(i)}|$ is uniformly bounded outside a finite set of indices. For $m>0$, $E\in\mathcal{ED}_m$ exactly when
\begin{itemize}
\item each $E\cap A_{(i)}$ belongs to $\mathcal{ED}_m\upharpoonright A_{(i)}$, and
\item the cardinalities $|E\cap A_{(i)}|$ are uniformly bounded outside finitely many indices.
\end{itemize}
Indeed, every generator other than a root selector is supported in one first-level block. Conversely, cover the finitely many exceptional intersections by their finitely many generators and the bounded remainder by finitely many root selectors.

Let $M_i=M\cap A_{(i)}$. There are infinitely many infinite $M_i$, and each of them is major below $(i)$. Pair each nonempty finite $M_i$ with a different infinite $M_j$; leave the remaining infinite intersections unpaired. This partitions $M$ into countably many infinite sets $D_j$, each consisting of one infinite $M_i$ and at most one finite $M_{i'}$.

For $m=0$, send these sets bijectively to the first-level blocks of $\omega$. The description above shows that this is an isomorphism: a uniform bound outside finitely many original blocks becomes at most twice that bound outside finitely many grouped blocks, and a bound on grouped blocks bounds each constituent block. Exceptional grouped blocks have at most two constituent blocks.

Now suppose $m>0$ and the assertion has been proved for $m-1$. Apply the induction hypothesis to the subsystem below $(i)$ and its major subset $M_i$. It gives an isomorphism between the local ideal restricted to $M_i$ and the local ideal on $A_{(i)}$. These are $\mathcal{ED}_m\upharpoonright M_i$ and $\mathcal{ED}_m\upharpoonright A_{(i)}$, respectively, by the preceding identification. Adjoining the paired finite set does not change its isomorphism type. To check this, take an infinite set in the local ideal, absorb the finitely many extra points into it by a bijection, and leave all other points fixed. Membership is unchanged because all changes occur inside a set belonging to the ideal.

Thus each $D_j$, with its induced ideal, can be mapped isomorphically onto a first-level block with the ideal $\mathcal{ED}_m\upharpoonright A_{(j)}$. Combine these bijections. The first condition in the membership description is preserved because it holds on a grouped block exactly when it holds on its infinite constituent. The second is preserved by the same factor-two bound used in the base case. This gives the required isomorphism. Its target is proper, so $M$ is positive.
\end{proof}

Consequently, if a dominant set is partitioned into finitely many pieces, one of the pieces contains a major $M$ on which each $\mathcal{ED}_m$ has the same isomorphism type as on $\omega$.

\begin{obs}\label{obs:one}
If $\mathcal I\approx\mathcal J$, on base sets $A$ and $B$ respectively, then for every $k\ge1$:
\begin{itemize}
\item $A\to(\mathcal I^+)^2_k$ if and only if $B\to(\mathcal J^+)^2_k$;
\item $A\to(<\omega,\ldots,<\omega\mid\mathcal I^+)^2_k$ if and only if $B\to(<\omega,\ldots,<\omega\mid\mathcal J^+)^2_k$.
\end{itemize}
\end{obs}

Both assertions follow by transporting a coloring along the isomorphism. We use them in the induction below.

\begin{lema}\label{lemma:three}
If $X\subseteq M$ is positive for $\mathcal{ED}_m\upharpoonright M$, where $M$ is major, then $X\in\mathcal{ED}_m^+$. Similarly, positivity for $\mathcal{ED}_{m+1}\upharpoonright A_{(n)}$ implies positivity for $\mathcal{ED}_{m+1}$.

Moreover, if $X\subseteq\omega$ satisfies
\[
(\forall r\in\omega)(\exists N\ge r)\quad |X\cap A_{(N)}|\ge r,
\]
then $X\in\mathcal{ED}_m^+$. The same conclusion holds if, for some $s\in\omega^{\le m}$, these inequalities hold with $A_{s^\frown N}$ in place of $A_{(N)}$.
\end{lema}
\begin{proof}
For a subset of the base of a restriction, membership in the restricted ideal is exactly membership in the original ideal. This proves the first two assertions.

For the last assertion, suppose that $X$ is covered by finitely many generators. At the root, all but the root selectors are supported in finitely many child blocks, and the root selectors give a uniform bound in every other child block. This contradicts the displayed condition. The same argument applies below $s$ when $|s|\le m$: terminal blocks and selectors below $s$ are confined to single child blocks, selectors at $s$ contribute at most one point per child, and selectors strictly above $s$ meet $A_s$ in at most one point. Thus outside finitely many children there is again a uniform bound.
\end{proof}

\begin{deff}
Following \cite[Definition 4.12]{articulomaicol}, we write $\omega\to(<\omega,\mathcal I^+)^2_2$ if every coloring $c:[\omega]^2\to2$ has either finite $0$-homogeneous sets of arbitrarily large size or an $\mathcal I$-positive $1$-homogeneous set.
\end{deff}

\begin{prop}\label{prop:one}
For an ideal $\mathcal I$ on $\omega$, the following are equivalent:
\begin{itemize}
\item $\omega\to(<\omega,\mathcal I^+)^2_2$;
\item $\omega\to(<\omega,<\omega\mid\mathcal I^+)^2_2$.
\end{itemize}
\end{prop}
\begin{proof}
Assume the first property and fix a coloring with no positive homogeneous set. Applying the property to the coloring and then to its complement gives finite homogeneous sets of arbitrarily large size in both colors. Conversely, apply the second property to a coloring. A positive $1$-homogeneous set gives the second alternative of the first property. A positive $0$-homogeneous set, or the finite-set alternative for both colors, gives arbitrarily large finite $0$-homogeneous sets.
\end{proof}

\begin{lema}\label{lemma:five}
If $X\in\mathcal{ED}_m^+$, there is a positive $Y\subseteq X$ such that $\mathcal{ED}_m\upharpoonright Y\approx\mathcal{ED}_{\mathrm{fin}}$.
\end{lema}
\begin{proof}
Choose a node $s$ of length at most $m$ such that $X'=X\cap A_s$ is positive but every $X'\cap A_{s^\frown i}$ belongs to $\mathcal{ED}_m$. To find it, start at the root and move to a positive child as long as one exists. The process stops before level $m+1$, since every block at that level is a generator.

For every $r$, the set $I_r=\{i:|X'\cap A_{s^\frown i}|\ge r\}$ is infinite. Otherwise, $X'$ would be covered by finitely many small child intersections and finitely many partial selectors at $s$, contradicting positivity. Choose distinct $k_n\in I_n$ for $n\ge1$, and choose $y_n\subseteq X'\cap A_{s^\frown k_n}$ of size $n$. Put $Y=\bigcup_{n\ge1}y_n$.

For $Z\subseteq Y$, we have
\[
Z\in\mathcal{ED}_m\quad\Longleftrightarrow\quad
\sup_{n\ge1}|Z\cap y_n|<\infty.
\]
For the forward direction, finitely many generators yield a uniform bound outside finitely many child blocks, as in Lemma~\ref{lemma:three}. The remaining intersections are finite because the $y_n$ are finite. For the reverse direction, a bound $r$ allows $Z$ to be covered by $r$ partial selectors at $s$. In particular $Y$ is positive. Any bijection sending $y_n$ to $\{(n,j):j<n\}$ is the required isomorphism with the ideal on $\Delta$.
\end{proof}

We now prove the positive partition relations for $\mathcal{ED}_m$. The next proposition is the main result of this section.

\begin{prop}\label{prop:ed-positive}
For every $m\in\omega$,
\[
\omega\to(\mathcal{ED}_m^+)^2_{m+1}
\quad\text{and}\quad
\omega\to(<\omega,\ldots,<\omega\mid\mathcal{ED}_m^+)^2_{m+2}.
\]
\end{prop}
\begin{proof}
We prove both assertions simultaneously by induction. For $m=0$, the first is the one-color relation. For the second, use $\mathcal{ED}^+\to(<\omega,\mathcal{ED}^+)^2_2$ from \cite[Theorem 4.14(4)]{articulomaicol}. Proposition~\ref{prop:ed-models}, Proposition~\ref{prop:one}, and Remark~\ref{obs:one} give the second assertion for $\mathcal{ED}_0$.

Suppose the assertions hold for $m=n$, and let $c:[\omega]^2\to n+2$. If $c$ has a positive homogeneous set for $\mathcal{ED}_{n+1}$, there is nothing to prove. Assume otherwise. For a major $M\subseteq A_{(i)}$ below $(i)$, apply Lemma~\ref{lemma:four} to the subsystem rooted at $(i)$, whose local ideal is $\mathcal{ED}_{n+1}\upharpoonright A_{(i)}$. Together with Lemma~\ref{lemma:two}, this gives $\mathcal{ED}_{n+1}\upharpoonright M\approx\mathcal{ED}_n$. On each such set the second inductive assertion therefore gives finite homogeneous sets of every color and every prescribed size: the positive alternative is excluded by our assumption on $c$.

\begin{lema}\label{lemma:seven}
Under this assumption on $c$, let $N\ge1$, let $D\subseteq A_{(k_0)}$ be major below $(k_0)$, and let $B\subseteq\omega\setminus\{k_0\}$ be infinite. For each $m\in B$, let $D_m\subseteq A_{(m)}$ be major below $(m)$. There are $a\in[D]^N$, an infinite $B'\subseteq B$, a color $i<n+2$, and major sets $D'_m\subseteq D_m$ below $(m)$ for $m\in B'$, such that $a$ is $i$-homogeneous and
\[
c(\{x,y\})=i\qquad
(x\in a,\ y\in\textstyle\bigcup_{m\in B'}D'_m).
\]
\end{lema}
\begin{proof}
Enumerate $D=\{x_j:j\in\omega\}$. At stage $j$, partition each of the remaining major sets by the color of its pairs with $x_j$. One piece is dominant below its first-level node, so take a major subset of that piece. A single color $i_j$ is obtained for infinitely many indices $m$; retain those indices as an infinite set $B_j$ and write $D_m^j$ for the retained major subsets. Make these choices nested: $B_{j+1}\subseteq B_j$ and $D_m^{j+1}\subseteq D_m^j$ when $m\in B_{j+1}$. Thus
\[
c(\{x_h,y\})=i_h\qquad
(h\le j,\ m\in B_j,\ y\in D_m^j).
\]
All pairs here have distinct points because $k_0\notin B$.

The sets $Z_i=\{x_j:i_j=i\}$ partition $D$. By Proposition~\ref{prop:local-dominant}, one of them contains a major set $D'$ below $(k_0)$. The observation preceding the lemma gives an $i$-homogeneous $a\in[D']^N$. Put $J=\max\{j:x_j\in a\}$, $B'=B_J$, and $D'_m=D_m^J$. The displayed invariant gives all the required cross-colors.
\end{proof}

Start with the major sets $A_{(m)}$ for all $m$. Repeatedly apply Lemma~\ref{lemma:seven}, at stage $j\ge1$ selecting an index $k_j$ from the remaining infinite set, discarding the indices at most $k_j$ from the reserve, and taking $N=j$. Retain the new major subsets supplied by the lemma for the next stage. This gives strictly increasing $k_j$, sets $a_j\subseteq A_{(k_j)}$ of size $j$, and colors $i_j$ such that $a_j$ is $i_j$-homogeneous and every pair from $a_j$ to a later $a_l$ has color $i_j$. Some color $i$ occurs for infinitely many $j$. Then $\bigcup_{i_j=i}a_j$ is $i$-homogeneous and positive by Lemma~\ref{lemma:three}, since its sizes are unbounded in distinct first-level blocks. This contradicts the assumption on $c$ and proves the first assertion for $m=n+1$.

For the second assertion, let $c:[\omega]^2\to n+3$. If every color has arbitrarily large finite homogeneous sets, we are done. Otherwise, rename as $0$ a color lacking that property, and merge colors $0$ and $1$ into a single color. The resulting coloring has $n+2$ colors, so the assertion just proved gives a positive homogeneous set $B$. If its color is not the merged color, it is already homogeneous for $c$.

Otherwise $c$ uses only colors $0$ and $1$ on $[B]^2$. By Lemma~\ref{lemma:five}, choose a positive $Y\subseteq B$ whose restricted ideal is isomorphic to $\mathcal{ED}_{\mathrm{fin}}$. Apply on $Y$ the asymmetric two-color relation for that ideal from \cite[Theorem 4.13(2)]{articulomaicol}. The alternative of arbitrarily large finite $0$-homogeneous sets is excluded by the choice of color $0$. We obtain a positive $1$-homogeneous set. This proves the second assertion for $m=n+1$ and completes the induction.
\end{proof}

\begin{teo}
For every $m\in\omega$, $\omega\to(\mathcal{ED}_m^+)^2_{m+1}$ but $\omega\not\to(\mathcal{ED}_m^+)^2_{m+2}$.
\end{teo}
\begin{proof}
Combine Propositions~\ref{prop:ed-negative} and \ref{prop:ed-positive}.
\end{proof}

Put $\mathcal{ED}_\omega=\bigcap_{m\in\omega}\mathcal{ED}_m$. It is a proper $F_{\sigma\delta}$ ideal, and it is tall because it contains $\mathcal{PC}$.

\begin{teo}
For every $n\ge1$, $\omega\to(\mathcal{ED}_\omega^+)^2_n$. For every function $f$ as in Definition~\ref{def:tc}, we also have $\omega\to(\mathcal{PC}(f)^+)^2_n$ and $\omega\to(\mathcal{TC}(f)^+)^2_n$.
\end{teo}
\begin{proof}
The ideal $\mathcal{ED}_\omega$ is contained in $\mathcal{ED}_{n-1}$. The same holds for $\mathcal{PC}(f)$ and $\mathcal{TC}(f)$: every tree in their respective additional generating families has finitely many nodes at level $n$, so its set is covered by finitely many blocks of that level. Apply Proposition~\ref{prop:ed-positive} to $\mathcal{ED}_{n-1}$. A positive set for that ideal is positive for all three smaller ideals.
\end{proof}

\begin{question}
Does $\omega\to(\mathcal{PC}^+)^n_k$ hold for every $n>2$ and every $k\ge2$?
\end{question}

The answer is negative for every $n\ge3$ and $k\ge2$. It is enough to construct a two-coloring of triples and then extend it to higher arities as below. Given three distinct points $a,b,c$, let $s$ be the longest node with $a,b,c\in A_s$. Either the three points lie in three distinct children of $s$, or two of them lie in one child and the third in another. Color the triple $0$ in the first case and $1$ in the second.

Let $H$ be $0$-homogeneous with at least two points. Fix $a,b\in H$ and let $s$ be the node at which they separate. For any other $c\in H$, the longest common node of $a,b,c$ is an initial segment of $s$; if it were shorter than $s$, then $a$ and $b$ would lie in one of its children and $c$ in another, and the triple would have color $1$. Hence $c$ lies in a third child of $s$. Thus $H$ is a partial selector at $s$ and belongs to $\mathcal{PC}$.

Let $H$ be $1$-homogeneous. No node has three children meeting $H$, so every node of $T_H$ has at most two children. Split $H$ according to the at most two first-level blocks it meets. Each part has one child at the root and at most two children at every deeper node, so it belongs to $\mathcal S_1$. Hence $H\in\mathcal{PC}$, and $\omega\not\to(\mathcal{PC}^+)^3_2$.

For $n>3$, color an $n$-element set by the color of its three least elements. An infinite homogeneous set for this coloring is homogeneous for the coloring of triples, since any triple in it can be completed by larger elements of the set; finite sets belong to $\mathcal{PC}$ anyway.

The same coloring works for $\mathcal{PC}(f)$. The $0$-homogeneous sets are again partial selectors. For a $1$-homogeneous $H$, choose $L$ with $f(L+1)\ge2$ and split $H$ into the finitely many blocks it meets at level $L$. Each part has at most one child at nodes of length less than $L$ and at most two children at deeper nodes, so it belongs to $\mathcal B(f)$. The coloring also works for every $\mathcal{ED}_m$ and for $\mathcal{ED}_\omega$, since a set meeting at most two children of every node is covered by finitely many blocks at any fixed level. It does not settle the case of $\mathcal{TC}(f)$ for slowly growing $f$: a set meeting at most two children of every node can have $2^n$ nodes at level $n$ and be $\mathcal{TC}(f)$-positive.

For other ideals satisfying the partition relations for pairs with every finite number of colors but failing the two-color relation for triples, see \cite[Propositions 5.13 and 5.16]{cano}. The coloring above gives this separation for $\mathcal{PC}$, $\mathcal{PC}(f)$, and $\mathcal{ED}_\omega$.

\section{Insights into the Solecki ideal}

We now study the Solecki ideal and a new ideal $\mathcal S_\omega$ extending it. Let $\mu$ be the usual product probability measure on $2^\omega$.

\begin{deff}
Let $\Omega=\{U\subseteq2^\omega:U\text{ is clopen and }\mu(U)=1/2\}$. The \emph{Solecki ideal} on the countable set $\Omega$ is
\[
\mathcal S=\{A\subseteq\Omega:\exists F\in[2^\omega]^{<\omega}\ \forall U\in A\ (U\cap F\ne\emptyset)\}.
\]
\end{deff}

The ideal $\mathcal S$ is tall and $F_\sigma$; see \cite{michaeldiagrama,solecki}. Its notation is unrelated to the selector families $\mathcal S_0,\mathcal S_1,\mathcal S_2$ used earlier. The examples in Section 1 answer the existence question for the cut-and-choose game, but they do not decide the corresponding question for $\mathcal S$, also posed by Zapletal.

\begin{question}
Does player II have a winning strategy in $G_1(\mathcal S)$?
\end{question}

We believe that the answer is negative, that is, that player I has a winning strategy in $G_1(\mathcal S)$, but we have no proof.

\begin{deff}
For $n\ge2$, let
\[
\Omega_n=\{V\subseteq2^\omega:V\text{ is clopen and }\mu(V)=2^{-n}\},
\]
and put
\[
\mathcal S_n^+=\{A\subseteq\Omega:\forall V\in\Omega_n\ \exists U\in A\ (U\cap V=\emptyset)\},
\qquad
\mathcal S_\omega^+=\bigcup_{n\ge2}\mathcal S_n^+.
\]
\end{deff}

A clopen set of measure at most $2^{-n}$ can be enlarged to one of measure exactly $2^{-n}$: express it as a union of cylinders at a sufficiently deep common level and add cylinders at that level. In particular, $\mathcal S_n^+\subseteq\mathcal S_{n+1}^+$. Every member of $\mathcal S_n^+$ is $\mathcal S$-positive: a finite set of points is contained in a clopen set of measure $2^{-n}$, which would meet every member of a family in $\mathcal S$. The superscript anticipates Proposition~\ref{prop:somega-ideal}, where these families become the positive sets of an ideal.

The definition also has the following useful interpretation. A family $A$ does not belong to $\mathcal S_\omega^+$ if and only if, for every $\varepsilon>0$, there is a clopen set $V$ with $\mu(V)<\varepsilon$ that intersects every member of $A$. In one direction choose $n$ with $2^{-n}<\varepsilon$; in the other enlarge a sufficiently small clopen set to a member of $\Omega_n$. We will use this interpretation to show that particular families are small.

\begin{lema}\label{lemma:six}
If $A\in\mathcal S_n^+$ and $A=B\cup C$, then $B\in\mathcal S_{n+1}^+$ or $C\in\mathcal S_{n+1}^+$.
\end{lema}
\begin{proof}
Otherwise choose $V_0,V_1\in\Omega_{n+1}$ intersecting every member of $B$ and $C$, respectively. Their union has measure at most $2^{-n}$. Enlarge it to a member of $\Omega_n$. The enlarged set intersects every member of $A$, contradicting $A\in\mathcal S_n^+$.
\end{proof}

\begin{prop}\label{prop:somega-ideal}
$\mathcal S_\omega=\mathcal P(\Omega)\setminus\mathcal S_\omega^+$ is an ideal.
\end{prop}
\begin{proof}
The family $\mathcal S_\omega^+$ is upward closed, and the preceding lemma shows that if a union is positive, one of its two parts is positive. Hence its complement is hereditary and closed under finite unions. Every finite family belongs to that complement: choose a point in each of its members and cover those finitely many points by a clopen set of arbitrarily small measure.

Finally, $\Omega\in\mathcal S_2^+$. Indeed, the complement of any clopen set of measure $1/4$ contains a clopen set of measure $1/2$. Thus $\mathcal S_\omega$ is proper.
\end{proof}

\begin{prop}\label{prop:somega-complexity}
$\mathcal S_\omega$ is a tall $F_{\sigma\delta}$ ideal containing $\mathcal S$.
\end{prop}
\begin{proof}
For a fixed clopen $V$, the family
\[
\{A\subseteq\Omega:\exists U\in A\ (U\cap V=\emptyset)\}
\]
is open in $2^\Omega$. Since $\Omega_n$ is countable, $\mathcal S_n^+$ is $G_\delta$. Its union over $n\ge2$ is $G_{\delta\sigma}$, so $\mathcal S_\omega$ is $F_{\sigma\delta}$.

If $A\in\mathcal S$ is witnessed by a finite $F\subseteq2^\omega$, a clopen neighborhood of $F$ of arbitrarily small measure intersects every $U\in A$. Thus $A\in\mathcal S_\omega$. Since $\mathcal S$ is tall, so is $\mathcal S_\omega$.
\end{proof}

\begin{deff}
Following \cite[Section 4.4]{articulomaicol}, for an ideal $\mathcal I$ on a countable set $X$, write $X\to(\omega,\mathcal I^+)^2_2$ if every $c:[X]^2\to2$ has an infinite $0$-homogeneous set or an $\mathcal I$-positive $1$-homogeneous set. Write $\mathcal I^+\to(\omega,\mathcal I^+)^2_2$ if
\[
A\to\bigl(\omega,(\mathcal I\upharpoonright A)^+\bigr)^2_2
\]
for every $A\in\mathcal I^+$.
\end{deff}

\begin{teo}\label{thm:somega-ramsey}
$\mathcal S_\omega^+\to(\omega,\mathcal S_\omega^+)^2_2$.
\end{teo}
\begin{proof}
First observe that if $A\in\mathcal S_n^+$ and $B\in\mathcal S_\omega$, then $A\setminus B\in\mathcal S_{n+1}^+$. Apply Lemma~\ref{lemma:six} to $A=(A\cap B)\cup(A\setminus B)$; the first part cannot belong to $\mathcal S_{n+1}^+$.

Fix $A\in\mathcal S_\omega^+$ and $c:[A]^2\to2$. For $D\subseteq A$ and $x\in D$, put $N_i^D(x)=\{y\in D\setminus\{x\}:c(\{x,y\})=i\}$. Start with $A_0=A$. If some $x_j\in A_j$ has $N_0^{A_j}(x_j)$ positive, choose such an $x_j$ and put $A_{j+1}=N_0^{A_j}(x_j)$. If this continues forever, the distinct points $x_j$ form an infinite $0$-homogeneous set.

Otherwise we reach a positive set $D$ such that $N_0^D(x)\in\mathcal S_\omega$ for every $x\in D$. Choose $m\ge2$ with $D\in\mathcal S_m^+$ and enumerate $\Omega_{m+1}=\{V_i:i\in\omega\}$. Recursively put
\[
B_i=D\setminus\left(\{y_j:j<i\}\cup\bigcup_{j<i}N_0^D(y_j)\right)
\]
and choose $y_i\in B_i$ with $y_i\cap V_i=\emptyset$. The removed family is in $\mathcal S_\omega$, so the initial observation gives $B_i\in\mathcal S_{m+1}^+$, making this choice possible. The points are distinct and pairwise of color $1$. Moreover, they form a member of $\mathcal S_{m+1}^+$, since for each $V_i$ the point $y_i$ is disjoint from it. This is the required positive homogeneous set.
\end{proof}

\begin{coro}
$\Omega\to(\omega,\mathcal S^+)^2_2$.
\end{coro}
\begin{proof}
Apply Theorem~\ref{thm:somega-ramsey} to $\Omega$ and use $\mathcal S\subseteq\mathcal S_\omega$.
\end{proof}

In particular, a coloring witnessing $\Omega\not\to(\mathcal S^+)^2_2$ must have infinite homogeneous sets of both colors. We next show that $\Omega\not\to(\mathcal S_\omega^+)^2_2$, and compare the two ideals using a related game.

For a related finite-selection game in which player I chooses an arbitrary positive set at each round, see \cite[Lemma 5.14]{survey}.

\begin{deff}\label{def:gfin}
The game $G_{\mathrm{fin}}(\mathcal I)$ is played as in $G_1(\mathcal I)$, except that player II chooses a finite subset of the retained piece instead of a point. More explicitly, player I first partitions $\omega=A_0^0\cup A_1^0$ and player II chooses $i_0\in2$ and $a_0\in[A_{i_0}^0]^{<\omega}$. At round $m\ge1$, player I partitions $A_{i_{m-1}}^{m-1}=A_0^m\cup A_1^m$, and player II chooses $i_m\in2$ and $a_m\in[A_{i_m}^m]^{<\omega}$. Player I wins if $\bigcup_m a_m\in\mathcal I$; otherwise player II wins. As in $G_1$, the retained piece must be nonempty, but empty finite selections are allowed.
\end{deff}

\begin{center}
\begin{tabular}{c|c|c|c|c}
I & $\omega=A_0^0\cup A_1^0$ & & $A_{i_0}^0=A_0^1\cup A_1^1$ & $\cdots$\\
II & & $i_0\in2,\ a_0\in[A_{i_0}^0]^{<\omega}$ & & $\cdots$
\end{tabular}
\end{center}

\begin{obs}
If player I wins $G_{\mathrm{fin}}(\mathcal I)$ and $\mathcal I\le_K\mathcal J$, then player I wins $G_{\mathrm{fin}}(\mathcal J)$. If player II wins $G_{\mathrm{fin}}(\mathcal I)$ and $\mathcal J\le_K\mathcal I$, then player II wins $G_{\mathrm{fin}}(\mathcal J)$.
\end{obs}

Use the simulations following Definition~\ref{def:game}, transporting finite sets in place of points.

\begin{prop}\label{prop:gfin-fsigma}
If $\mathcal I$ extends to an $F_\sigma$ ideal, then player II has a winning strategy in $G_{\mathrm{fin}}(\mathcal I)$.
\end{prop}
\begin{proof}
Let $\mathcal I\subseteq\mathcal J=\Fin(\varphi)$, where $\varphi$ is a lower semicontinuous submeasure. At round $n$, player II keeps a $\mathcal J$-positive part of the current $\mathcal J$-positive piece. This part has infinite submeasure, so lower semicontinuity provides a finite $a_n$ in it with $\varphi(a_n)>n$. The union of the $a_n$ has infinite submeasure and is not in $\mathcal J$, hence not in $\mathcal I$.
\end{proof}

\begin{deff}
The ideal $\mathrm{conv}$ on $\mathbb Q\cap[0,1]$ consists of the sets with only finitely many accumulation points in $[0,1]$. It is a tall $F_{\sigma\delta\sigma}$ ideal \cite{michaeldiagrama}.
\end{deff}

\begin{prop}
Player I has a winning strategy in $G_{\mathrm{fin}}(\mathrm{conv})$.
\end{prop}
\begin{proof}
Maintain a retained piece $D_n$ contained in a closed interval $J_n$ of length $2^{-n}$, starting with $D_0=\mathbb Q\cap[0,1]$ and $J_0=[0,1]$. If $t_n$ is the midpoint of $J_n$, partition $D_n$ into its points at most $t_n$ and its points greater than $t_n$. Let $D_{n+1}$ be the part chosen by player II, and let $J_{n+1}$ be the corresponding closed half of $J_n$. This is a legal partition even when $D_n$ omits an endpoint. The intervals $J_n$ have a unique common point $x$.

For every neighborhood of $x$, all sufficiently late finite selections lie in that neighborhood. The finitely many earlier selections contain only finitely many points. Thus their union has no accumulation point other than possibly $x$, and belongs to $\mathrm{conv}$.
\end{proof}

\begin{deff}\label{def:binary-choice}
For $s\in2^{<\omega}$, write $[s]=\{x\in2^\omega:s\subseteq x\}$. For $U\in\Omega$, define $x_U\in2^\omega$ recursively. Having chosen $s=x_U\upharpoonright n$, put
\[
x_U(n)=
\begin{cases}
0,&\mu(U\cap[s^\frown0])\ge\mu(U\cap[s^\frown1]),\\
1,&\text{otherwise}.
\end{cases}
\]
\end{deff}

Every chosen cylinder has positive intersection with $U$. Since $U$ is closed, $x_U\in U$. Since it is also open, some cylinder along $x_U$ is contained in $U$; from that stage onward the two halves have equal measure and the rule always chooses $0$. Thus $x_U$ is eventually zero.

\begin{prop}\label{prop:two}
Player I wins $G_{\mathrm{fin}}(\mathcal S_\omega)$, whereas player II wins $G_{\mathrm{fin}}(\mathcal S)$.
\end{prop}
\begin{proof}
The assertion for player II follows from Proposition~\ref{prop:gfin-fsigma}, since $\mathcal S$ is $F_\sigma$.

For player I, put $D_s=\{U\in\Omega:s\subseteq x_U\}$. Starting from $D_\emptyset=\Omega$, partition $D_s$ into $D_{s^\frown0}$ and $D_{s^\frown1}$. The choices of player II determine a branch $x\in2^\omega$ and finite selections $a_n\subseteq D_{x\upharpoonright(n+1)}$.

Fix $\varepsilon>0$ and choose $L\ge1$ with $2^{-L}<\varepsilon/2$. For $n\ge L$, every $U\in a_n$ contains its point $x_U\in[x\upharpoonright L]$, so that cylinder intersects all members of the tail. There are only finitely many members in $\bigcup_{n<L}a_n$. Choose a small cylinder inside each of them, with the sum of their measures less than $\varepsilon/2$. The union of those cylinders and $[x\upharpoonright L]$ is clopen, has measure less than $\varepsilon$, and intersects every member of $\bigcup_n a_n$. Since this works for every $\varepsilon$, the selected family belongs to $\mathcal S_\omega$.
\end{proof}

In particular, the inclusion $\mathcal S\subseteq\mathcal S_\omega$ is strict.

\begin{prop}\label{prop:conv-somega}
$\mathrm{conv}\le_K\mathcal S_\omega$.
\end{prop}
\begin{proof}
Let $\pi:2^\omega\to[0,1]$ be the continuous binary map
\[
\pi(x)=\sum_{n\in\omega}x(n)2^{-n-1},
\]
and put $f(U)=\pi(x_U)$. Each $x_U$ is eventually zero, so $f$ maps $\Omega$ into $\mathbb Q\cap[0,1]$.

Take $B\in\mathrm{conv}$. Its closure in $[0,1]$ is countable and compact: it consists of $B$ together with finitely many accumulation points. Since every real has at most two binary expansions, $C=\pi^{-1}[\overline B]$ is also countable and compact. Hence $\mu(C)=0$. For every $\varepsilon>0$, cover $C$ by an open set of measure less than $\varepsilon$ and then, by compactness, by a finite union $V$ of cylinders inside that open set. For each $U\in f^{-1}[B]$, the point $x_U$ belongs to $C\cap U\subseteq V\cap U$. Thus $V$ intersects every such $U$, and the small-clopen characterization gives $f^{-1}[B]\in\mathcal S_\omega$.
\end{proof}

\begin{coro}
$\Omega\not\to(\mathcal S_\omega^+)^2_2$.
\end{coro}
\begin{proof}
The reduction $\mathcal R\le_K\mathrm{conv}$ \cite[Remark 4.3]{articulomaicol}, together with Proposition~\ref{prop:conv-somega}, gives $\mathcal R\le_K\mathcal S_\omega$. Apply Proposition~\ref{prop:second}.
\end{proof}

Together with Theorem~\ref{thm:somega-ramsey}, this shows that $\mathcal S_\omega$ is a tall $F_{\sigma\delta}$ ideal satisfying $\mathcal S_\omega^+\to(\omega,\mathcal S_\omega^+)^2_2$ but $\mathcal S_\omega^+\not\to(\mathcal S_\omega^+)^2_2$.

\begin{deff}
The ideal $\mathrm{nwd}$ on $\mathbb Q$ consists of the subsets that are nowhere dense in $\mathbb Q$. It is a tall $F_{\sigma\delta}$ ideal \cite{balzar,michaeldiagrama}.
\end{deff}

\begin{deff}
In the ideal game $G_3(\mathcal I)$ (see \cite{laflame} and \cite[proof of Theorem 3.4.1]{tesisdavid}), for an ideal on a countable set $X$, player I chooses $I_k\in\mathcal I$ at round $k$, and player II chooses $n_k\in X\setminus I_k$. Player I wins if $\{n_k:k\in\omega\}\in\mathcal I$; otherwise player II wins.
\end{deff}

\begin{prop}[M. Hru\v{s}\'{a}k and D. Meza-Alc\'{a}ntara, {\cite[proof of Theorem 3.4.1]{tesisdavid}}]\label{prop:g3-nwd}
Let $\mathcal I$ be a Borel ideal. If player II has a winning strategy in $G_3(\mathcal I\upharpoonright X)$ for every $X\in\mathcal I^+$, then $\mathcal I\le_K\mathrm{nwd}$.
\end{prop}

\begin{prop}
$\mathcal S_\omega\le_K\mathrm{nwd}$.
\end{prop}
\begin{proof}
Let $X\in\mathcal S_\omega^+$ and choose $n\ge2$ with $X\in\mathcal S_n^+$. Enumerate $\Omega_{n+1}=\{V_k:k\in\omega\}$. At round $k$ of $G_3(\mathcal S_\omega\upharpoonright X)$, player I chooses $I_k\in\mathcal S_\omega\upharpoonright X$. As in the first observation of Theorem~\ref{thm:somega-ramsey}, Lemma~\ref{lemma:six} gives $X\setminus I_k\in\mathcal S_{n+1}^+$. Choose $U_k\in X\setminus I_k$ disjoint from $V_k$. The selected family $\{U_k:k\in\omega\}$ belongs to $\mathcal S_{n+1}^+$, so this strategy wins. Proposition~\ref{prop:g3-nwd} now applies to every positive $X$.
\end{proof}

We have obtained $\mathcal S\subsetneq\mathcal S_\omega$ and
\[
\mathcal R\le_K\mathrm{conv}\le_K\mathcal S_\omega\le_K\mathrm{nwd},
\]
while $\mathcal S_\omega$ still satisfies the asymmetric relation in Theorem~\ref{thm:somega-ramsey}. The question for $\mathcal S$ itself asks whether $\mathcal R\le_K\mathcal S$, equivalently whether its two-color Ramsey property fails. We have worked on this question at length without success; in our opinion it needs new ideas.

\section{Remarks on \texorpdfstring{$K$}{K}-Uniform Ideals}

An ideal $\mathcal I$ is \emph{$K$-uniform} if $\mathcal I\upharpoonright X\le_K\mathcal I$ for every $X\in\mathcal I^+$. Recall that, up to isomorphism, $\mathcal{ED}_{\mathrm{fin}}$ is the restriction of $\mathcal{ED}$ to $\Delta=\{(i,j):j<i\}$. The following properties are known; see \cite{michaeldiagrama,articulomaicol}:
\begin{itemize}
\item $\mathcal{ED}_{\mathrm{fin}}$ is tall and $F_\sigma$;
\item $\mathcal{ED}_{\mathrm{fin}}$ is $K$-uniform;
\item if $\mathcal I$ is a tall $F_\sigma$ $K$-uniform ideal, then $\mathcal{ED}_{\mathrm{fin}}\le_K\mathcal I$: the proof of \cite[Corollary 4.4]{articulomaicol} gives an $\mathcal I$-positive $X$ with $\mathcal{ED}_{\mathrm{fin}}\le_K\mathcal I\upharpoonright X$, and $K$-uniformity completes the reduction.
\end{itemize}
As far as we know, $\mathcal{ED}_{\mathrm{fin}}$ was the only known tall $F_\sigma$ $K$-uniform ideal. We construct another one, not Kat\v{e}tov-equivalent to $\mathcal{ED}_{\mathrm{fin}}$. We identify a graph on $\omega$ with its edge set, a subset of $[\omega]^2$; the edge set may be empty.

\begin{deff}
A proper vertex coloring of a graph $G\subseteq[\omega]^2$ is a map $c:\omega\to\kappa$ such that $c(a)\ne c(b)$ whenever $\{a,b\}\in G$. Its \emph{chromatic number} $\chi(G)$ is the least number of colors in such a coloring.
\end{deff}

\begin{deff}[See \cite{tesisdavid,michaeldiagrama}]
The ideal $\mathcal G_{\mathrm{fc}}$ on $[\omega]^2$ consists of graphs of finite chromatic number. Equivalently, it is generated by the bipartite graphs.
\end{deff}

This is a tall $F_\sigma$ ideal, and $\mathcal S\le_K\mathcal G_{\mathrm{fc}}$ and $\mathcal{ED}_{\mathrm{fin}}\le_K\mathcal G_{\mathrm{fc}}$; see \cite[Theorem 1.6.22]{tesisdavid}.

\begin{deff}
Let $\langle I_n:n\ge1\rangle$ partition $\omega$, with $|I_n|=n$. Let $K_n=[I_n]^2$ and $K=\bigcup_{n\ge1}K_n$. Define $\mathcal K$ on $K$ by
\[
\mathcal K=\{X\subseteq K:\text{for some }r\ge2,\ X\text{ contains no clique on }r\text{ vertices}\}.
\]
\end{deff}

This is an ideal. Heredity and membership of finite sets are immediate. For finite unions, if $X$ has no clique of size $r$ and $Y$ has none of size $s$, a clique of size $R(r,s)$ in $X\cup Y$ would, by the finite Ramsey theorem, contain a clique of size $r$ in $X$ or one of size $s$ in $Y$. To apply that theorem, color an edge by $0$ when it belongs to $X$, and by $1$ otherwise. Finally, $K\notin\mathcal K$, since it contains cliques of every finite size.

\begin{obs}
$\mathcal G_{\mathrm{fc}}\le_K\mathcal G_{\mathrm{fc}}\upharpoonright K$.
\end{obs}

The inclusion $K\hookrightarrow[\omega]^2$ witnesses this reduction. The restriction is proper because $\chi(K)$ is infinite.

\begin{teo}\label{thm:k-uniform}
$\mathcal K$ is a tall $F_\sigma$ $K$-uniform ideal that is not Kat\v{e}tov-equivalent to $\mathcal{ED}_{\mathrm{fin}}$.
\end{teo}
\begin{proof}
Let $A\in\mathcal K^+$. We choose copies $K'_n$ of $K_n$ in $A$, for $n\ge2$, lying in distinct original blocks. At each stage only finitely many original blocks have been used. Their union is finite, so removing them from $A$ leaves cliques of arbitrarily large size. Hence a new block contains a clique of the required size.

For each $n\ge2$, choose a graph isomorphism from $K_n$ to $K'_n$ and combine the induced maps on edges into $h:K\to A$. If $B\in\mathcal K\upharpoonright A$ has no clique on $r$ vertices, neither does $h^{-1}[B]$: a clique in the preimage lies within one original block, where $h$ is induced by a vertex injection. Thus $h^{-1}[B]\in\mathcal K$, proving $\mathcal K\upharpoonright A\le_K\mathcal K$. This establishes $K$-uniformity.

For tallness, an infinite set of edges meets infinitely many of the finite blocks $K_n$. Choose one edge from each of infinitely many such blocks. The resulting infinite graph has no clique of size $3$ and belongs to $\mathcal K$.

For each $r\ge2$, the family of subsets of $K$ containing no clique on $r$ vertices is closed. Indeed, a clique witnessing nonmembership is a finite set of $\binom r2$ edges, and every graph containing those edges remains outside the family. The union of these closed families over $r$ is $\mathcal K$, so $\mathcal K$ is $F_\sigma$.

Finally, every graph of finite chromatic number has bounded clique size. Therefore
\[
\mathcal S\le_K\mathcal G_{\mathrm{fc}}
\le_K\mathcal G_{\mathrm{fc}}\upharpoonright K
\le_K\mathcal K.
\]
Also, $\mathcal S\nleq_K\mathcal{ED}_{\mathrm{fin}}$. Indeed,
\[
\mathcal{ED}_{\mathrm{fin}}\le_K\mathcal I_{1/n}\subseteq\mathcal Z,
\]
where $\mathcal I_{1/n}$ is the summable ideal and $\mathcal Z$ is the density-zero ideal. The latter has the Fubini property, whereas every ideal Kat\v{e}tov above $\mathcal S$ fails that property; see \cite[Section 2]{michaeldiagrama}. Thus $\mathcal S\le_K\mathcal{ED}_{\mathrm{fin}}$ would imply $\mathcal S\le_K\mathcal Z$, a contradiction. Consequently $\mathcal K\nleq_K\mathcal{ED}_{\mathrm{fin}}$, and the two ideals are not Kat\v{e}tov-equivalent.
\end{proof}

This answers Hru\v{s}\'{a}k's question. Once the ideal $\mathcal K$ is written down, the proof is short; the work was in finding the right object. We expect that the same construction with $n$-uniform hypergraphs in place of graphs gives, for different $n$, tall $F_\sigma$ $K$-uniform ideals that are pairwise Kat\v{e}tov-inequivalent, but we have not checked this.

\paragraph{Acknowledgments.}
The author thanks Michael Hru\v{s}\'{a}k for his supervision and guidance. This work forms part of the author's PhD thesis, written between 2014 and 2017.

\paragraph{Funding.}
This research was supported by the Mexican National Council for Science and Technology (CONACYT).

\end{document}